\documentclass[11pt, reqno]{amsproc}

\usepackage{latexsym}
\usepackage{amssymb} 
\usepackage{mathrsfs}
\usepackage{amsmath}
\usepackage{latexsym}
\usepackage{delarray}
\usepackage{amssymb,amsmath,amsfonts,amsthm,
mathrsfs}

\RequirePackage{srcltx}

\usepackage{hyperref}
\renewcommand\eqref[1]{(\ref{#1})} 

\theoremstyle{plain}
\newtheorem{thm}{Theorem}

\newtheorem{cor}[thm]{Corollary}

\newtheorem{remark}[thm]{Remark}

\newcommand{\mb}[1]{\ensuremath{\mathbb{#1}}}

\newcommand{\R}{\mb{R}}

\newcommand{\T}{\mb{T}}

\def\Dcal{{\mathcal D}}

\newcommand{\dd}{ {\rm d}}
\newcommand{\Gh}{ {\widehat{G}}}
\newcommand{\Lap}{\mathcal{L}}

\def\dxi{{d_\xi}}

\def\HS{{\mathtt{HS}}}
\def\whf{\widehat{f}}
\def\whfhs{\|\widehat{f}(\xi)\|_{\HS}}

\newcommand{\jp}[1]{\langle{#1}\rangle}
\newcommand{\p}[1]{\left({#1}\right)}

\newcommand{\Tr}{{\mathrm{Tr}}}
 \newcommand{\supp}{\operatorname{supp}}

 \newcommand{\conv}{\operatorname{conv}}

\title[Nikolskii inequality and functional classes on compact Lie groups
]{
Nikolskii inequality and functional classes on compact Lie groups
}

\author[E. D. Nursultanov, M. V. Ruzhansky, S. Yu. Tikhonov]{E. D. Nursultanov, M. V. Ruzhansky, S. Yu. Tikhonov}
\address{E. D. Nursultanov:
\endgraf
Moscow State University (Kazakhstan branch)
\endgraf
and Eurasian National University
\endgraf
Astana
\endgraf
Kazakhstan
\endgraf
{\it E-mail address} {\rm er-nurs@yandex.ru}
\endgraf
\,\qquad
\endgraf
M. V. Ruzhansky:
 \endgraf
  Department of Mathematics
  \endgraf
  Imperial College London
  \endgraf
  180 Queen's Gate, London SW7 2AZ
  \endgraf
  United Kingdom
  \endgraf
  {\it E-mail address} {\rm m.ruzhansky@imperial.ac.uk}
  \endgraf
\,\qquad
\endgraf
S. Yu. Tikhonov:
 \endgraf
ICREA
and Centre de Recerca Matem\`{a}tica (CRM)  \endgraf
 E-08193, Bellaterra  \endgraf
 Barcelona
  \endgraf
  {\it E-mail address} {\rm stikhonov@crm.cat}
}

\thanks{E. Nursultanov was supported by grants MOH PK 4080 GF-4,  3311 GF-4.
M. Ruzhansky was supported by EPSRC Leadership Fellowship, EPSRC грант
EP/K039407/1 and Leverhulme Grant RPG-2014-02. S. Tikhonov was supported by
grants 2014-SGR-289 from AGAUR (Generalitat de
Catalunya), MTM 2014-59174-P, and RFFI-13-01-00043.
}

\date{}

\begin{document}


\begin{abstract}
In this note we study  Besov, Triebel--Lizorkin, Wiener, and Beurling 
function spaces on compact Lie groups.
A major role in the analysis is played by the Nikolskii inequality. 
\end{abstract}

\maketitle

{\bf{1. Introduction.}} \ \  The classical Nikolskii inequality for trigonometric polynomials $T_L$ of degree up to $L$ can be written as (\cite{nikol}):
$$
 \|T_L\|_{L^q(\T)} \le 2 L^{1/p-1/q} \|T_L\|_{L^p(\T)},
$$
where $1\le p<q\le \infty$.
 The Nikolskii inequality plays an important role in the analysis of different function spaces (for example, see \cite{triebel}) and in
 the approximation theory (for example, see  \cite{ditzian}).

In the Euclidean case, for functions $ f\in L^p(\R^n)$ such that $\supp (\widehat{f})$ is compact (cf. \cite{nessel}) we have
\begin{equation}\label{nik-first}
 \|f\|_{L^q(\R^n)} \le \left(C(p) \mu(\conv [\supp (\widehat{f})]) \right)^{1/p-1/q}  \|f\|_{L^p(\R^n)},
\end{equation}
where $1\le p\le q\le \infty$, $\mu(E)$ is the Lebesgue measure of $E$, and
$\conv [E]$ is the convex hull of $E$. Inequalities of the form (\ref{nik-first}) are often called
Plancherel--Polya--Nikolskii inequalities.

Recently, in \cite{pes} Pesenson obtained the Bernstein--Nikolskii inequality on symmetric spaces of noncompact type, and in \cite{Pesenson:Besov-2008} on compact homogeneous spaces.

Let $G$  be a compact Lie group of dimension $\dim G$ and let $\Gh$ be its unitary dual.
If we fix bases in representation spaces we can work with matrix representations
$\xi:G\to\mb{C}^{d_{\xi}\times d_{\xi}}$ of dimensions $d_{\xi}$.
By the Peter--Weyl theorem the system
$\{\sqrt{d_{\xi}}\xi_{ij}: [\xi]\in\Gh, 1\leq i,j\leq d_{\xi}\}$
is an orthonormal basis in
 $L^{2}(G)$ with respect to the normalized Haar measure on $G$.
All the integrals below and the spaces $L^{p}(G)$ will be always considered with respect to this
normalized bi-invariant Haar measure on $G$.

For  $f\in C^{\infty}(G)$ we define its Fourier coefficient at $\xi\in [\xi]\in\Gh$ by
$$
\widehat{f}(\xi)=\int_{G} f(x) \xi(x)^{*} \dd  x.
$$
Thus, we have  $\widehat{f}(\xi)\in\mb{C}^{d_{\xi}\times d_{\xi}}.$
The Fourier series of a function $f$ takes the form
\begin{equation}\label{EQ:FS}
f(x)=\sum_{[\xi]\in\Gh} d_{\xi} \Tr (\widehat{f}(\xi) \xi(x)).
\end{equation}

For  $[\xi]\in\Gh$ by $\jp{\xi}$ we denote the eigenvalue of the operator
$(I-\Lap_{G})^{1/2}$ corresponding to the representation class  $[\xi]\in\Gh$, where $\Lap_{G}$ is the Laplacian on $G$, see,
for example \cite[Chapter 1.7]{stein}.

In \cite{rt:book} the following Lebesgue spaces $\ell^{p}(\Gh)$ on $\Gh$ were defined as follows:
using the Fourier coefficients of $f$, we set
\begin{equation}\label{EQ:norm}
\|\widehat{f}\|_{\ell^{p}(\Gh)}=\left(\sum_{[\xi]\in\Gh} d_{\xi}^{p(\frac{2}{p}-\frac12)}
\|\widehat{f}(\xi)\|_{\HS}^{p}\right)^{1/p},\; 1\leq p<\infty,
\end{equation}
and
\begin{equation}\label{EQ:norm-linfty}
\|\widehat{f}\|_{\ell^{\infty}(\Gh)}=
\sup_{[\xi]\in\Gh} d_{\xi}^{-\frac12} \|\widehat{f}(\xi)\|_{\HS},
\end{equation}
 where
$\|\widehat{f}(\xi)\|_{\HS}=\Tr(\widehat{f}(\xi) \widehat{f}(\xi)^{*})^{1/2}.$
For these spaces the following Hausdorff--Young inequalities are valid:
$$
\|\widehat{f}\|_{\ell^{p'}(\Gh)}\leq \|f\|_{L^{p}(G )},\;
\|f\|_{L^{p'}(G )}\leq \|\widehat{f}\|_{\ell^{p}(\Gh)},
\quad
1\leq p\leq 2,\quad \frac{1}{p'}+\frac1p=1.
$$

Let $N(L)$ be the Weyl eigenvalue counting function for the elliptic
pseudo-differential operator
$(1-\Lap_{G })^{1/2}$, denoting the number of its eigenvalues
$\leq L$ counted with multiplicities. Then
$$
N(L)=\mathop{\sum_{\jp{\xi}\leq L}}_{[\xi]\in\Gh} d^2_{\xi}.
$$
For sufficiently large $L$ the Weyl asymptotic formula says that
$$
N(L)\sim C_{0} L^{n}, \quad C_{0}=(2\pi)^{-n}\int_{\sigma_{1}(x,\omega)< 1} \dd x\dd\omega,
$$
where $n=\dim G $, and the integral is taken with respect to the canonical measure on the cotangent bundle
 $T^{*}(G )$ induced by the canonical symplectic form. Here
 $\sigma_{1}$ is the principal symbol of the operator $(1-\Lap_{G })^{1/2},$ see e.g. \cite{shubin:r}.

\medskip
Full proofs of our results below will appear in \cite{nrt15}.

\medskip

{\bf{2. Nikolskii inequality.}} \ \ Let $T$ be a trigonometric polynomial on a compact Lie group $G$, i. e. a function with
only finitely many non-zero Fourier coefficients. Let $D$ be the Dirichlet kernel, i.e. the function
 $D\in C^{\infty}(G)$ such that
$$
\widehat{D}(\xi):= I_{d_{\xi}}\quad \textrm{ for } \; \jp{\xi}\leq L,
$$
and zero otherwise. Here $I_{d_{\xi}}\in \mathbb{C}^{{d_{\xi}}\times {d_{\xi}}}$ denote the identity matrix.

\begin{thm}\label{THM3}
Let $0<p<q\leq\infty$. For $0<p\leq 2$ set $\rho:=1$, and for $2<p<\infty$ let $\rho$ be the smallest integer $\geq p/2$.
Then
$$
\|T\|_{L^q(G)}\leq \p{\sum_{\widehat{T^\rho}(\xi)\not=0} d_\xi^2}^{\frac 1p-\frac 1q} \|T\|_{L^p(G)}.
$$
Moreover, this inequality is sharp for $p=2$ and $q=\infty$, and the equality is attained at
$T=D$.
\end{thm}
We note that for the classical trigonometric polynomials of several variables the Nikolskii inequality is
well known (\cite{nikol}).

\begin{remark}
Note that if $T=T_{L}$, i.e. if $\widehat{T}(\xi)=0$ for $\jp{\xi}>L$,
then $\sum_{\widehat{T}(\xi)\not=0} d_\xi^2\leq N(L)$ and, therefore,
$$
\|T_{L}\|_{L^{q}(G)}\leq N(\rho L)^{\frac1p-\frac1q}\|T_{L}\|_{L^{p}(G)}
\asymp
\big({\rho}L\big)^{n(\frac1p-\frac1q)}\|T_{L}\|_{L^{p}(G)}.
$$
\end{remark}
For a partial sum of the Fourier series of $f$:
$$
S_{L}f(x)=\sum_{\jp{\xi}\leq L} d_{\xi}\ \Tr(\whf(\xi)\xi(x))
$$
one can prove the following result.

\begin{cor}\label{THM3-2}
Let  $G$ be a compact Lie group and let
$1\leq p<q\leq\infty$ be such that $\frac 1p > \frac 1q + \frac12.$
Then we have
$$
\left(\sum_{k=1}^\infty\frac{\left(k^{1-1/p+1/q}\sup_{N(L)\geq k}\frac1{N(L)}\|S_{L}f\|_{L^q(G)}\right)^p}k\right)^{1/p}\leq C \|f\|_{L^p(G)}
$$
for all $f\in L^{p}(G)$.
In particular, we have
$\ \  {N(L)}^{\frac 1q-\frac 1p}\|S_{L}f\|_{L^q(G)} =o(1)$ as $L\to \infty.$
\end{cor}

\medskip

{\bf{3. Embeddings of functional classes.}} \ \
Here we investigate embedding theorems and interpolation properties of several classes of functions on
a compact Lie group $G$.
Using the definition \eqref{EQ:FS} of the Fourier series, we can defined Sobolev, Besov, and Triebel--Lizorkin spaces,
respectively,
as follows:
\begin{equation*}\label{EQ:Sobolev}
H^{r}_{p}=H^{r}_{p}(G )=\left\{f\in\Dcal'(G ):\,\, \|f\|_{H^{r}_{p}} :=
  \Big\|(1-\Lap_{G })^{r/2}f\Big\|_p<\infty\right\},
\end{equation*}
\begin{equation*}\label{EQ:Besov1}
B^{r}_{p,q}=B^{r}_{p,q}(G )=\left\{f\in\Dcal'(G ):\,\, \|f\|_{B^{r}_{p,q}} :=
 \left(\sum_{s=0}^\infty 2^{sr q} \Big\|\sum_{2^{s}\le\jp{\xi}<2^{s+1}}
 \dxi\ \Tr\p{\whf(\xi)\xi(x)}\Big\|_p^q\right)^{1/q}<\infty\right\},
\end{equation*}
$$
F^{r}_{p,q}=F^{r}_{p,q}(G )=\left\{f\in\Dcal'(G ):\,\,  \|f\|_{F^{r}_{p,q}} :=
\Bigg\| \Bigg(\sum_{s=0}^\infty 2^{sr q} \, \Big|\sum_{2^{s}\le\jp{\xi}<2^{s+1}}
 \dxi\ \Tr\p{\whf(\xi)\xi(x)}\Big|^q\Bigg)^{1/q}\Bigg\|_p<\infty\right\}.
$$

Then we have the following result:

\begin{thm}\label{THM:Besov}
Let $G $ be a compact Lie group of dimension $n$.
Then
\begin{itemize}

\item[(1)] \ \
$B_{p_1,q}^{r_1}\hookrightarrow B_{p_2,q}^{r_2}, \qquad
0<p_1\le p_2\leq\infty,\; 0<q\le \infty,\; r_2= r_1 - n(\frac{1}{p_1}-\frac{1}{p_2});$

\item[(2)] \ \
$B_{p,\min\{p, 2\}}^{r}\hookrightarrow H^r_p \hookrightarrow B_{p,\max\{p, 2\}}^{r}, \qquad r\in\R,\quad 1< p<\infty;$

\item[(3)] \ \
$B_{p,q}^r\hookrightarrow L_q, \qquad
1<p<q<\infty, \; r=n(\frac1p-\frac1q)$;

\item[(4)] \ \
$B_{p,1}^r\hookrightarrow L_\infty, \qquad
0<p\le\infty,\quad
 \; r=\frac np$;

\item[(5)] \ \
$B_{p,\min\left\{p,q\right\}}^{r}\hookrightarrow F_{p,q}^{r}\hookrightarrow B_{p,\max\left\{p,q\right\}}^{r},
\qquad
1<p<\infty,\;
0<p<\infty,\; 0<q\le\infty$.

\item[(6)] \ \
$(B_{p,\beta_0}^{r_0},B_{p,\beta_1}^{r_1})_{\theta ,q} = (H_{p}^{r_0},H_{p}^{r_1})_{\theta ,q} = (F_{p,\beta_0}^{r_0},F_{p,\beta_1}^{r_1})_{\theta ,q} = B_{p,q}^r,
\qquad
0<r_1 < r_0 < \infty$, $0< \beta_0, \beta_1, q \le \infty$,
$1<p<\infty,$ $ r=(1-\theta)r_0+\theta r_1, \qquad 0<\theta<1.$

\end{itemize}
\end{thm}

For functions on the torus the corresponding results can be found, for example, in the book \cite{triebel}.

Consequently, using norms (\ref{EQ:norm}) and (\ref{EQ:norm-linfty}), we can investigate the embeddings
between Wiener and Beurling classes defined as follows:
$$
A^\beta(\Gh )= 
\left\{f\in\Dcal'(G ):\,\,
\|f\|_{A^\beta}:=\|\whf\|_{\ell^{\beta}(\Gh)}=
\p{\sum_{[\xi]\in\Gh} d_{\xi}^{\beta(\frac2\beta-\frac12)}\whfhs^{\beta}}^{1/\beta}
<\infty\right\}
$$
and
$$
A^{*,\beta}(\Gh)=\left\{f:\,\,\|f\|_{A^{*,\beta}(\Gh)}:=\left(\sum_{s=0}^\infty 2^{ns} \left(\sup\limits_{2^s\le \jp{\xi}}
d_{\xi}^{-1/2}\whfhs\right)^\beta\right)^{1/\beta}<\infty\right\},
$$
where  $0<\beta<\infty$. For periodic functions, i.e. for
$G= \mathbb{T}^n$, we have $d_{\xi}\equiv 1$,
$\Gh\simeq \mathbb{Z}^n$, and
 $\whfhs = |{\widehat{f}(\xi)}|$,
and such spaces have been investigated, for example, in  \cite{beu, lifl} and \cite[Ch. 6]{trigub}.

\begin{thm}\label{THM:Wiener-beta}
Let $G $ be a compact Lie group of dimension $n$.
\\
{\textnormal {(A)}}.\ \
Let $\alpha>0$ and $\frac 1\beta=\frac n\alpha+\frac{1}{p'}$.
Then
\begin{equation}\label{sz}
\|f \|_{A^\beta} \le C  \|f\| _{B^{\alpha }_{p,\beta}},\qquad 1< p\le 2;
\end{equation}
\[
 \|f\| _{B^{\alpha}_{p,\beta}}  \le C \|f \|_{A^\beta}, \qquad 2\le p <\infty.
\]
{\textnormal {(B)}}.\ \
Let $0<\beta<\infty$ и $p\geq 2$.
Then
$$
C_1 \|f\|_{B^{n(\frac1\beta-\frac{1}{p'})}_{p,\beta}}
  \le
\|f\|_{A^{*,\beta}}\le C_2 \|f\|_{B^{n/\beta}_{1,\beta}}
$$
\end{thm}

As a consequence, we obtain
$$
C_1\|f\|_{A^{\beta}} \le
 \|f\|_{B^{n(1/\beta-1/2)}_{2,\beta}}
  \le  C_2
\|f\|_{A^{*,\beta}}.
$$
The left inequality is an analogue of Bernstein theorem on the absolute convergence of Fourier series.
It strengthens the following inequality proved by Faraut in \cite{faraut} for groups $G$ of unitary matrices:
if $f\in C^{k}(G)$ for an even $k>\frac{\dim G}{2}$ then $\whf\in\ell^{1}(\Gh)$, i.e. $f\in A(G)$.
For periodic functions of several variables inequality (\ref{sz}) 
follows from the results in \cite{szasz}.

Finally, we look at the following Beurling-type spaces:
$$
A_r^{*,\beta}=\left\{f: \|f\|_{A_r^{*,\beta}}
:=\left(\sum_{s=0}^\infty\Big(2^{r ns}\sup\limits_{2^s\le \jp{\xi}}
d_{\xi}^{-1/2}\|\widehat{f}(\xi)\|_{\HS}\Big)^\beta\right)^{\frac1\beta}<\infty\right\}.
$$
Note that $A_{1/\beta}^{*,\beta}=A^{*,\beta}$. These spaces are interpolation spaces in the following sense:

\begin{thm} \label{THM:interpolation}
Let $0<r_1 < r_0 < \infty$, $0< \beta_0, \beta_1, q \le \infty$, $ r=(1-\theta)r_0+\theta r_1$ и $0<\theta<1.$
Then
$$
(A_{r_0}^{*,\beta_0},
A_{r_1}^{*,\beta_1}
)_{\theta, q} = A_r^{*,q}.
$$
In particular,
$$
(A^{*,1/r_0},A^{*,1/r_1})_{\theta, 1/r} = A^{*,1/r}.
$$
\end{thm}

\bibliographystyle{alphaabbr}

\bibliography{bib-Nikolski}

\end{document}